\documentclass{amsart}
\usepackage{amssymb,amsfonts,amscd}
\numberwithin{equation}{section}
\let\cal\mathcal
\def\Ascr{{\cal A}}
\def\Bscr{{\cal B}}

\def\Fscr{{\cal F}}

\def\Lscr{{\cal L}}

\def\Oscr{{\cal O}}

\let\blb\mathbb
\def\CC{{\blb C}} 
 
\def\FF{{\blb F}} 
\def\QQ{{\blb Q}}
\def\GG{{\blb G}}

\def \ZZ{{\blb Z}}

\def \NN{{\blb N}}
\def \RR{{\blb R}}
\def \HH{{\blb H}}

\def\Stab{\operatorname{Stab}}
\def\kar{\operatorname{char}}
\def\re{\operatorname{Re}}

\def\pr{\mathop{\text{pr}}\nolimits}

\def\quot{/\!\!/}

\def\mod{\operatorname{mod}}

\def\gr{\operatorname{gr}}
\def\Lie{\operatorname{Lie}}

\def\Supp{\mathop{\text{\upshape Supp}}}

\def\rad{\operatorname {rad}}

\def\gr{\operatorname {gr}}
\def\Spec{\operatorname {Spec}}

\def\Gl{\operatorname {Gl}}
\def\Rep{\operatorname {Rep}}

\def\Ext{\operatorname {Ext}}
\def\Hom{\operatorname {Hom}}
\def\End{\operatorname {End}}

\def\Tr{\operatorname {Tr}}

\def\End{\operatorname {End}}
\def\Gal{\operatorname {Gal}}

\def\rk{\operatorname {rk}}

\def\r{\rightarrow}

\def\Irr{\operatorname{Irr}}
\def\udim{\underline{\operatorname{dim}}}
\def\choose#1#2{\begin{pmatrix} #1\\ #2 \end{pmatrix}}

\newtheorem{lemma}{Lemma}[section]
\newtheorem{proposition}[lemma]{Proposition}
\newtheorem{theorem}[lemma]{Theorem}

\newtheorem{lemmas}{Lemma}[subsection]
\newtheorem{propositions}[lemmas]{Proposition}
\newtheorem{theorems}[lemmas]{Theorem}
\newtheorem{corollarys}[lemmas]{Corollary}

\newtheorem*{theoremstar}{Theorem}

\theoremstyle{definition}

\newtheorem{convention}[lemma]{Convention}
\newtheorem{conjecture}[lemma]{Conjecture}

\newtheorem{definitions}[lemmas]{Definition}

{

}

\theoremstyle{remark}

\newdimen\uboxsep \uboxsep=1ex
\def\uboxn#1{\vtop to 0pt{\hrule height 0pt depth 0pt\vskip\uboxsep
\hbox to 0pt{\hss #1\hss}\vss}}

\def\uboxs#1{\vbox to 0pt{\vss\hbox to 0pt{\hss #1\hss}
\vskip\uboxsep\hrule height 0pt depth 0pt}}

\title[Absolutely indecomposable representations]{Absolutely
indecomposable representations  and Kac-Moody Lie algebras
(with an appendix by Hiraku Nakajima)}
\author{William Crawley-Boevey}
\address{
Department of Pure Mathematics\\
University of Leeds\\
Leeds LS2\, 9JT\\
UK} 
\email{w.crawley-boevey@leeds.ac.uk}
\author{Michel Van den Bergh}
 \address{Departement WNI,  Limburgs Universitair Centrum, 
 3590 Diepenbeek, Belgium.}
  \email{vdbergh@luc.ac.be}
\thanks{The second author is a senior researcher at the FWO}
\dedicatory{Dedicated to Idun Reiten on the occasion of her sixtieth birthday.}
\begin{document}
\begin{abstract} 
 A conjecture of Kac states that the polynomial counting the number of absolutely
  indecomposable representations of a quiver over a finite field with
  given dimension vector has
  positive coefficients and furthermore that its constant term is
  equal to the multiplicity of the corresponding root in the
  associated Kac-Moody Lie algebra. In this paper we prove these
  conjectures for indivisible dimension vectors.
 \end{abstract}
\maketitle

\section{Introduction}
Let $Q$ be a finite quiver without loops with vertices  $I$ and fix  $\alpha\in\NN^I$. In
\cite{Kac1} V. Kac showed (over an algebraically closed field) that $Q$ has an
indecomposable representation of dimension vector $\alpha$ if and
only if $\alpha$ is a root of a certain Kac-Moody Lie algebra $\mathfrak{g}$
associated to $Q$. This was a spectacular generalization of  earlier
results by Gabriel \cite{Gabriel1} for the finite type case and Dlab and
Ringel \cite{DlabRingel} for the tame case.

Now assume that the ground field is finite. In this case one should
consider \emph{absolutely indecomposable representions}, i.e.
indecomposable representations which remain indecomposable over
the algebraic closure of the ground field.

For  $\alpha\in\NN^I$ let $a_\alpha(q)$
be the number of absolutely indecomposable representations of $Q$
with dimension vector $\alpha$
over $\FF_{q}$. Kac has shown that $a_\alpha(q)$ is a
polynomial in $q$ with integral coefficients \cite{kac2}.
Regarding this polynomial Kac made the following  intriguing
conjectures:
{\def\thelemma{A}
\begin{conjecture} $a_\alpha(q)\in \NN[q]$.
\end{conjecture}
}
{\def\thelemma{B}
\begin{conjecture} If $\alpha$ is a root then $a_\alpha(0)$ is the
  multiplicity of $\alpha$ in $\mathfrak{g}$.
\end{conjecture}
} Despite our greatly increased understanding of the relationship
between quivers and Kac-Moody Lie algebras (thanks to Ringel, Lusztig,
Kashiwara, Nakajima and others) and despite the fact that over twenty years have
passed since these conjectures were stated, virtually no progress has
been made towards their proof. See \cite{hua1,Lebruyn:quivers,sev} for some partial
and related
results. 

In this paper we make the first substantial progress by proving the
following result:
\setcounter{lemma}{0}
\begin{theorem}
Conjecture A and B
  are true if $\alpha$ is indivisible.
\end{theorem}
To prove such a  result it is clear that 
one should first find a good cohomological interpretation for the
polynomial $a_\alpha(q)$. Unfortunately the equivariant cohomology of the
representation space of $Q$ (which is the obvious choice) counts
representations with multiplicity (see \cite{Behrend, Kim}) and
this yields trivial results in our case. 

Thus one of the main results in this paper is a new interpretation of
$a_\alpha(q)$ in the case that $\alpha$ is indivisible.
To state
this new interpretation we have to introduce some notations. We assume
temporarily that our base field is $\CC$. 
The
double $\bar{Q}$ of $Q$ is the quiver obtained by adding a
reverse arrow $a^\ast:j\r i$ for each arrow $a:i\r j$ in $Q$. The
preprojective algebra of $Q$ is $\Pi^0=\CC\bar{Q}/(\sum [a,a^\ast])$
where the sum runs over the arrows in $Q$. 

Define a bilinear form on $\CC I$ by $i\cdot j=\delta_{ij}$ and let
$\lambda\in\ZZ I$ be such that $\lambda\cdot \alpha=0$ but
$\lambda\cdot \beta\neq 0$ for $0<\beta<\alpha$. Then we show in
\S\ref{ref-2-4} that
\begin{equation}
\label{ref-1.1-2}
a_\alpha(q)=\sum_{i=0}^d \dim H^{2d-2i}(X_s,\CC)\,q^i
\end{equation}
where $X_s$ is the (smooth) moduli-space of $\lambda$-stable
$\Pi^0$-representations of dimension vector $\alpha$ \cite{King} and $d=1/2 \dim
X_s$. It is clear that this formula proves Conjecture A for
indivisible $\alpha$.

Now let
$\Lambda_\alpha=\Rep(\Pi^0,\alpha)^{\text{nil}}$  be the nilpotent
representations in the representation space of $\alpha$-dimensional
representations of $\Pi^0$. Lusztig has shown  \cite[Thm
12.9]{Lusztig2}\cite{Lusztig5} that $\Lambda_\alpha$ is a 
Lagrangian subvariety of the affine space $\Rep(\bar{Q},\alpha)$ and
furthermore that the irreducible components of $\Lambda_\alpha$ index a basis of
$U(\mathfrak{g}^+)_\alpha$ (see also \cite{SK}). We first observe that
Conjecture B for $\alpha$ indivisible is equivalent to the following.
\begin{proposition}
\label{ref-1.5-3}
 Let $\alpha$ be indivisible. The number of irreducible
  components of $\Lambda_\alpha$ which contain a $\lambda$-stable (or
  equivalently: semistable)
  representation is equal to $\dim \mathfrak{g}_\alpha$.
\end{proposition}
We then prove this proposition by relating the Harder-Narasimhan filtration
on $\Pi^0$-representations to the PBW-theorem for
$U(\frak{g}^+)$. This approach was partially suggested by a talk of
M.\ Reineke. See \cite{Reineke}.

Let us now sketch how we prove \eqref{ref-1.1-2}. Unless otherwise specified
our base field is now finite.
We show first
that  $a_\alpha(q)$ counts the points of a smooth
affine variety $X$ related to a \emph{deformed} preprojective algebra of $Q$
\cite{CH}. Our aim is then to count the points on $X$ using the Lefschetz fixed
point formula for the Frobenius action on $l$-adic cohomology.

 Since we are not able to
extract any meaningful results directly from  $X$, we
construct a one-parameter family $\Xi$ of smooth varieties whose general fiber is $X$ and
whose special fiber is $X_s$. Now it is easy to see $X_s$ carries a $\GG_m$-action whose fixed point set
is projective. By combining the Weil conjectures with results
from \cite{Bia,Bia1} we deduce from this that 
the absolute values of the eigenvalues of the Frobenius action on the cohomology of $X_s$ are the same
as those of
a smooth projective variety (see Appendix \ref{ref-A-23}).

 Since $\Xi$
is not locally trivial we cannot directly transfer results from $X_s$
to $X$. However  an argument involving the  hyper-K\"ahler
structure on the representation space of  $\bar{Q}$ 
shows  that $X_s$ and $X$ are homeomorphic for the analytic topology
in characteristic zero (see \cite[Cor.\ 4.2]{Nakajima}). By specialization this implies that $X_s$ and
$X$ have  isomorphic cohomology in large
characteristic. Unfortunately it is not immediately clear to us that this
isomorphism is compatible with Frobenius (think of the example given
by elliptic curves).

Therefore we refine Nakajima's  argument in such a way that it shows
that the family $\Xi$ is
trivial for the analytic topology  (see
lemma \ref{ref-2.3.3-15} below).
It follows that the cohomology of the fibers of $\Xi$ is constant in
large characteristic. Thus $X$ and $X_s$ have the same
cohomology even when the  Frobenius action is taken into account.
 This allows us to prove \eqref{ref-1.1-2} using
a simple technical lemma (see lemma \ref{ref-A.1-24}).

Some words on the organization of this paper. The proof of
\eqref{ref-1.1-2} and the equivalence of Conjecture B and
Proposition \ref{ref-1.5-3} are contained in Section \ref{ref-2-4}.  
The proof of \eqref{ref-1.1-2}
relies on a few basic results on
$l$-adic cohomology and invariant theory over $\ZZ$. We have collected
those in two appendices so that they don't detract from the main
arguments. The proof of Proposition 1.2 is contained in Section
\ref{newsection}.

We wish to thank Henning Andersen for some useful information
regarding invariants over $\ZZ$. We also wish to thank Markus Reineke
for communicating us the main results of \cite{Reineke}.

At the end of the paper we include an appendix by H. Nakajima which 
avoids the arguments of Section~\ref{cohomsection} by showing directly 
that two varieties have the same number of points over finite fields. 
We have retained the original Section~\ref{cohomsection}, however, 
since it shows more---the existence of a canonical isomorphism between 
the cohomology of 
$\operatorname {Rep}(\Pi^\lambda,\alpha)^\lambda/\!\!/ G(\alpha)$ 
and $\operatorname {Rep}(\Pi^0,\alpha)^\lambda/\!\!/ G(\alpha)$ 
for arbitrary $\lambda$ and~$\alpha$ (see below for notations).

\section{Proof of \eqref{ref-1.1-2} and the equivalence of Proposition
  \ref{ref-1.5-3} and conjecture B}
\label{ref-2-4}
\subsection{Notations and constructions}
\label{ref-2.1-5}
Let $Q=(I,Q,h,t)$  be a finite quiver without loops with vertices $I$ and edges
$Q$. $h,t$ are the maps which associate starting and ending vertex to an
edge. There is a standard symmetric bilinear form on $\ZZ^I$ given by
\[
(i,j)=
\begin{cases}
1&\text{if $i=j$}\\
-\frac{1}{2}\#\{\text{arrows between $i$ and $j$}\} & \text{if $i\neq j$}
\end{cases}
\]
We let $\mathfrak{g}$ be the Kac-Moody Lie algebra whose Cartan matrix
$(a_{ij})_{ij}$ is given by $a_{ij}=2(i,j)$. 

An absolutely indecomposable representation
 of $Q$ over a field $k$ is an indecomposable representation $V$ with the
property that $V\otimes_k \bar{k}$ is indecomposable, or equivalently $\End(V)/\rad \End(V)=k$.
For $\alpha\in\NN^I$, $a_\alpha(q)$
is  the number isomorphism classes of absolutely indecomposable
 representations of $Q$ with dimension vector $\alpha$
over the finite field $\FF_{q}$.

We now introduce some standard constructions related to the quiver
$Q$. Since we want to  use lifting to characteristic zero we need to
define things over $\ZZ$. This makes our notations a little pedantic
for which we apologize in advance. For some basic material with
respect to invariants over $\ZZ$ we refer to Appendix \ref{ref-B-31}.
The essential ingredient, \emph{on which we will rely tacitly below}, is that
all constructions are compatible with base change over an open part of
$\Spec \ZZ$.

Let $\bar{Q}$ be the double quiver of $Q$. Thus $\bar{Q}$ has the same
vertices as $Q$ but the edges are given by $\{a,a^\ast\mid a\in Q\}$
where $h(a^\ast)=t(a)$ and $t(a^\ast)=h(a)$. 

If $R$ is a commutative ring and $\lambda\in R^I$ then $\Pi^\lambda$ is the corresponding deformed
preprojective algebra \cite{CH}. Thus
\begin{equation}
\label{ref-2.1-6}
\Pi^\lambda=R\bar{Q}/\left(\sum_{a\in Q}[a,a^\ast]-\sum_{i\in I}\lambda_i i\right)
\end{equation}
For $\alpha,\beta\in\NN$ let $M_{\alpha\times\beta}$, $M_\alpha$,
$\Gl(\alpha)$ be  the $\ZZ$-schemes corresponding respectively to the
$\alpha\times\beta$-matrices, the $\alpha\times\alpha$-matrices and
the invertible $\alpha\times\alpha$-matrices.

For $\alpha\in \NN^I$ we define
$
\Rep(Q,\alpha)=\prod_{e\in Q} M_{\alpha_{h(e)}\times \alpha_{t(e)}}
$.
 We use  corresponding notations for $\bar{Q}$ and
$\Pi^\lambda$. 

For $i,j\in I$ put 
$i\cdot j=\delta_{ij}$. This defines a bilinear form on $R^I$ for any
ring $R$.
\begin{lemmas}
\label{ref-2.1.1-7}
If $R$ is a field  and 
if $\alpha\cdot \lambda\neq 0$ in $R$ then
$\Rep(\Pi^\lambda,\alpha)=\emptyset$.
\end{lemmas}
\begin{proof}
This follows from the standard trace argument.
\end{proof}
We also define 
$
\Gl(\alpha)=\prod_{i\in I} \Gl(\alpha_i)
$
and
we put $G(\alpha)=\Gl(\alpha)/\GG_m$. 

The Lie algebra of $\Gl(\alpha)$ is given by $M(\alpha)=\prod_i
M_{\alpha_i\times \alpha_i}$.  Over a field $l$ 
 we may
identify $\Lie(\Gl(\alpha)_l)$ with its dual via the trace
pairing. Under this pairing the 
dual to $\Lie(G(\alpha)_l)$ is identified with the trace zero matrices
in $M(\alpha)_l$. We denote the variety of trace zero matrices with $M(\alpha)^0$.

The algebraic group $G(\alpha)$ acts by conjugation on
$\Rep(Q,\alpha)$ and the orbits $\Rep(Q,\alpha)(l)/G(\alpha)(l)$ for
$l$ a field correspond
to isomorphism classes of $Q$-representa\-tions defined over $l$.

Now let $\lambda\in\ZZ^I$ such that $\lambda\cdot\alpha=0$. Then
$\lambda$ defines a character $\chi_\lambda$ of $G(\alpha)$ given by $(x_i)_{i\in
  I}\mapsto \prod_i\det(x_i)^\lambda_i$. As in \cite{King}, $\chi$
defines a line bundle $\Lscr$ on $\Rep(\bar{Q},\alpha)$. 
 We define
$\Rep(\bar{Q},\alpha)^{\lambda}$ as the $\Lscr$-semistable part \cite[\S
II]{seshadri} of $\Rep(\bar{Q},\alpha)$. Using the Hilbert-Mumford
criterion \cite[Prop. 3.1]{King} one finds that if $k$ is an
algebraically closed field then $V\in \Rep(\bar{Q},\alpha)(k)$ lies in
$\Rep(\bar{Q},\alpha)(k)^\lambda$ if and only if
\begin{equation}
\label{ref-2.2-8}
\lambda\cdot \underline{\text{\upshape dim}}V'\ge 0
\end{equation}
for every subrepresentation $0\neq V'\subsetneq V$. If we replace the
 inequality in \eqref{ref-2.2-8} by a strict one then we obtain the
 stable representations.  

Consider the  map 
\begin{equation}
\label{ref-2.3-9}
\mu: \Rep(\bar{Q},\alpha)\r M(\alpha)^0:(x_a)_{a\in \bar{Q}}\mapsto \sum
[x_a,x_a^\ast]_{a\in Q} 
\end{equation}
Over a field $l$, $\mu$ may be identified with a suitable moment map
for the $G(\alpha)_l$ action on $\Rep(\bar{Q},\alpha)_l$ via the identification
of $\Lie(G(\alpha))_l^\ast$ with
$M(\alpha)_l^0$.  We will refer to \eqref{ref-2.3-9} as \emph{the}
moment map. We clearly have $\mu^{-1}(\lambda)=\Rep(\Pi^\lambda,\alpha)$.

Let $L$ be the line in the affine space in $M(\alpha)^0$ spanned by
$0$ and $\lambda$ and let $W=\mu^{-1}(L)\cap
\Rep(\bar{Q},\alpha)^\lambda$.  Put $\Xi=W\quot G(\alpha)$ and let $f:\Xi\r L$ be the
induced map.
We put $X=f^{-1}(\lambda)=\Rep(\Pi^\lambda,\alpha)^\lambda\quot G(\alpha)$ and
$X_s=f^{-1}(0)=\Rep(\Pi^0,\alpha)^\lambda\quot G(\alpha)$.
\begin{definitions} We say that $\lambda\in\ZZ^I$ is \emph{generic}
  with respect to $\alpha\in \NN^I$ if  $\lambda\cdot \alpha=0$ but
$\lambda\cdot
  \beta\neq 0$ for all $0<\beta<\alpha$ (note that such a 
 $\lambda$ exists if and only if  $\alpha$ is indivisible).
\end{definitions}
If $\lambda$ is generic for $\alpha$ then it follows from
\eqref{ref-2.2-8} that over an algebraically closed 
field the notions of $\lambda$-semi-stability and $\lambda$-stability
coincide.
\begin{lemmas}
\label{ref-2.1.3-10}
Assume that $\lambda$ is generic with respect to $\alpha$.
Then there exists a non-empty open $U\subset \Spec \ZZ$ such that  $\Rep(\Pi^\lambda,\alpha)_U^\lambda=\Rep(\Pi^\lambda,\alpha)_U$.
\end{lemmas}
\begin{proof}
It is sufficient to prove this over $k=\bar{\QQ}$. In that case every
$x\in \Rep(\Pi^\lambda,\alpha)(k)$ is simple by lemma
\ref{ref-2.1.1-7}. Then by  \eqref{ref-2.2-8}
 it follows that $x$ is semistable (in fact stable) for
$\lambda$. 
\end{proof}
Since we are only interested in large characteristics we will commit a
slight abuse of notation by identifying
$X$ with $\Rep(\Pi^\lambda,\alpha)\quot G(\alpha)$ in the case that
$\lambda$ is generic. This is justified by the last lemma. 
\begin{lemmas} Assume that $\lambda$  is generic with respect to $\alpha$. Then there exists a non-empty open $U\subset \Spec \ZZ$ such that  the map $f:\Xi_U\r L_U$ is smooth.
\end{lemmas}
\begin{proof}
Again it is sufficient to do this over $k=\bar{\QQ}$. 

First we note that if $x\in \Rep(\bar{Q},\alpha)^\lambda(k)$ then by
\eqref{ref-2.2-8} 
$\End(x)=k$ and in particular $G(\alpha)_k$ acts freely on
$\Rep(\bar{Q},\alpha)^\lambda_k$. 

By  lemma \ref{ref-2.1.5-11} below $\mu$ is smooth
at $x$. Thus the restriction of $\mu$ to
$\Rep(\bar{Q},\alpha)^\lambda_k$ is smooth. It follows that the
induced map $W_k\r L_k$ is
also smooth. 

Since $G(\alpha)_k$ acts freely on $W_k$ we deduce that $W_k\r
W_k/G(\alpha)_k=\Xi_k$ is also smooth. This then yields that $\Xi_k\r
L_k$ is surjective on tangent spaces and hence smooth.
\end{proof}
We have used the following standard lemma. 
\begin{lemmas} 
\label{ref-2.1.5-11}
Let $X$ be a smooth symplectic variety over an
  algebraically closed field $k$ 
  and
  assume that $G$ is a linear algebraic group acting symplectically on
  $X$. Assume that in addition there is a moment map $X\mapsto
  \mathfrak{g}$ where $\mathfrak{g}=\Lie(G)$. Let $x\in X$. If the
  differential in $x$ of the $G$-action 
  $\mathfrak{g}\r T_x(X)$ is injective then $\mu$ is smooth at $x$.
\end{lemmas}

\subsection{Reformulation of Kac's conjectures for indivisible
  dimension vectors}

We assume throughout that $\alpha\in\NN^I$ is indivisible.
We put $k=\bar{\FF}_p$ and we let $q$ be a power of  $p$. We
prove the following result.
\begin{propositions}
\label{ref-2.2.1-12}
  Assume that $\lambda\in \ZZ^I$  is generic for $\alpha\in \NN^I$ and let $X=\Rep(\Pi^\lambda,\alpha)\quot G(\alpha)$ be as in \S\ref{ref-2.1-5}. Then for $p\gg 0$
we have
\[
a_\alpha(q)=q^{-d} |X(\FF_q)|
\]
with $d=1-(\alpha,\alpha)$
\end{propositions}
\begin{proof}
We consider the projection map
\[
\pi:\Rep(\Pi^\lambda,\alpha)\r \Rep(Q,\alpha)
\]
According to  \cite[Thm 3.3]{CB} the image of $\pi(\FF_q)$ consists
of indecomposable representations.
Since $\alpha$ is indivisible, representations of dimension vector
$\alpha$ are absolutely indecomposable if and only if they are
indecomposable.
Thus the image of $\pi(\FF_q)$ consists of absolutely
indecomposable representations.  

Let $\Rep(Q,\alpha)^{a.i}$ denote the constructible subset of 
absolutely indecomposable representations in the affine space $\Rep(Q,\alpha)$.
It is also shown in 
loc. cit. that the elements of $\Rep(Q,\alpha)^{a.i.}(\FF_q)$ lift to
$\Rep(\Pi^\lambda,\alpha)$. More precisely the inverse image of  $x\in\Rep(Q,\alpha)^{a.i.}(\FF_q)$ 
can be identified with $\Ext^1(x,x)^\ast$.

Starting from a variant of the Burnside formula we compute
\begin{align*}
\left|\Rep(Q,\alpha)^{a.i.}(\FF_q)/G(\alpha)(\FF_q)\right|&=
\frac{1}{|G(\alpha)(\FF_q)|} \sum_{x\in \Rep(Q,\alpha)^{a.i}(\FF_q)} |\Stab_{G(\alpha)}(x)|\\
&=q^{-1}\frac{1}{|G(\alpha)(\FF_q)|} \sum_{x\in \Rep(Q,\alpha)^{a.i}(\FF_q)}
|\End(x)|\\
&=q^{-1}\frac{1}{|G(\alpha)(\FF_q)|} \sum_{x\in \Rep(\Pi^\lambda,\alpha)(\FF_q)}
\frac{|\End(x)|}{|\Ext^1(x,x)|}\\
&=q^{(\alpha,\alpha)-1} \frac{|\Rep(\Pi^\lambda,\alpha)(\FF_q)|}{|G(\alpha)(\FF_q)|}
\end{align*}
 where we have used that $(-,-)$ is the
 symmetrization of the Euler form on $K_0(\mod(kQ))$.

Since $p\gg 0$ the inequalities defining genericity
also hold in $\FF_p$.  Hence we will assume this. 
By lemma \ref{ref-2.1.1-7} our choice of $\lambda$ insures that $\Rep(\Pi^\lambda,\alpha)(k)$
contains only simple representations. Thus if $x\in
\Rep(\Pi^\lambda,\alpha)(\FF_q)$ then $\End(x)=\FF_q$ and hence
$x$ has trivial stabilizer in $G(\alpha)(k)$.

Using
\cite[Cor. 5.3.b]{KR} we obtain
\begin{align*}
|\Rep(\Pi^\lambda,\alpha)(\FF_q)|/|G(\alpha)(\FF_q)|&=
|\Rep(\Pi^\lambda,\alpha)(\FF_q)/G(\alpha)(\FF_q)|\\ &=
|(\Rep(\Pi^\lambda,\alpha)(k)/G(\alpha)(k))^{\Gal(k/\FF_q)}|\\
&=|X(k)^{\Gal(k/\FF_q)}|=|X(\FF_q)| \qed
\end{align*}
\def\qed{}\end{proof}

\subsection{Cohomological triviality}
\label{cohomsection}
According to the program outlined in the introduction we want to
compare the cohomology of $X$ and $X_s$ (see \S\ref{ref-2.1-5}). One way to do this
is to show that $R^if_{!}(\QQ_l)$ is constant, at least over an open
part of the base  $\Spec\ZZ$.   This is the content of the next
proposition. Note that we do not assume that $\lambda$ is generic with
respect to $\alpha$.
\begin{propositions} 
\label{ref-2.3.1-13}
 There exists a non empty open $U\subset \Spec \ZZ$
  such that for every~$i$, $R^i{f_!}(\QQ_l)_U$ is the pullback of a sheaf
  on $U$.
\end{propositions}
\begin{corollarys}
\label{ref-2.3.2-14}
Let $k=\overline{\FF}_p$. For $p\gg 0$ there is an isomorphism between
$H^i_c(X_{s,k},\QQ_l)$ and $H^i_c(X_k,\QQ_l)$ which is compatible with the Frobenius action.
\end{corollarys}
\begin{proof} Let $f_s$, $f_g$ be the restrictions of $f$ to $X_s$
  and $X$.

Using the previous proposition and the fact that $R^if_!$
  commutes with base change we find for $p\gg 0$:
  $R^if_{s!,\FF_p}(\QQ_l)\cong R^if_{g!,\FF_p}(\QQ_l)$ on $\Spec \FF_p$. We may consider
  $R^if_{s!,\FF_p}(\QQ_l)$ and $R^if_{g!,\FF_p}(\QQ_l)$ as the
  $\Gal(k/\FF_p)$-modules given by 
  $H^i_c(X_{s,k},\QQ_l)$ and $H^i_c(X_k,\QQ_l)$ respectively. Since the
  Frobenius action is determined by the action of $\Gal(k/\FF_p)$
  \cite[\S 1.8]{Rapport} this proves
  what we want.
\end{proof}

\begin{proof}[Proof of Proposition \ref{ref-2.3.1-13}]   We use Deligne's generic base change result for direct images
  \cite[Thm 1.9]{Finitude}. This result was only stated for torsion
  sheaves, but the corresponding result for $l$-adic sheaves is an easy
  consequence. 

Since $f$ is of finite type there are only a finite
  number of $i$ for which $R^i{f_!}(\QQ_l)$ is non-zero. So we may
  treat each $i$ separately. Put $\Fscr= R^i{f_!}(\QQ_l)$. Let $g:L\r
  \ZZ$ be the structure map and let $\epsilon: g^\ast g_\ast
  \Fscr\r \Fscr$ be the map given by adjointness. Let $\Ascr,\Bscr$ be the 
  kernel and cokernel of $\epsilon$. By \cite[Thm 1.9]{Finitude} $g^\ast g_\ast
  \Fscr$ and hence $\Ascr,\Bscr$ will be constructible over an open subset 
  $V\subset \Spec \ZZ$.

Below we show  that $\epsilon_\CC: g_\CC^\ast g_{\CC,\ast} \Fscr_\CC\r \Fscr_\CC$
is an isomorphism. By \cite[Thm 1.9]{Finitude} we have $g_\CC^\ast
g_{\CC,\ast} \Fscr_\CC=(g^\ast g_\ast \Fscr)_\CC$. Hence
$\Ascr_\CC=\Bscr_\CC=0$. {}From the fact that  $\Ascr_V$ and $\Bscr_V$ are constructible
it follows
that $\Supp (\Ascr_V)$ and $\Supp(\Bscr_V)$ are constructible subsets
of $\Xi$ whose image in $\Spec \ZZ$ does not contain the generic point.
Hence we find $\Ascr_U=\Bscr_U=0$ for a suitable  open $U\subset V$.

Now we prove our claim that $\epsilon_\CC$ is an isomorphism.
To do this we replace the etale topology on $\Xi_\CC, L_\CC$ with the analytic
topology.
 Then the claim follows from the comparison theorem  \cite[\S
6.1.2]{BBD}, lemma \ref{ref-2.3.3-15} below and the fact that $L_\CC$ is connected.
\end{proof}
\emph{In the rest of this subsection our base field will be $\CC$ so we
  drop the corresponding subscript.}
\begin{lemmas}
\label{ref-2.3.3-15} 
$f:\Xi\r L$ is a trivial (topological) family.
\end{lemmas}
\begin{proof}
Let $V=\Rep(\bar{Q},\alpha)$. We will use the hyper-K\"ahler structure
on $V$ which was introduced by Kronheimer \cite{Kronheimer}. For the benefit of the
reader we recall the basic facts. First we define a Riemannian metric on $V$ via the
trace form:
\begin{equation}
\label{ref-2.4-16}
(x,y)=\re \sum_{a\in \bar{Q}}\Tr(x_a y_a^\dagger)
\end{equation}
where $z^\dagger$ is the conjugate transpose to $z$.

Let $\HH=\RR+\RR I+\RR J +\RR K$ be the
quaternions. We define an action of $\HH$ on $V$ via
\begin{align*}
I(x_a)_{a\in
  \bar{Q}}&=(ix_a)_{a\in\bar{Q}}\\
 J(x_a,x_{a^\ast})_{a\in
  Q}&=(-x^\dagger_{a^\ast},x^\dagger_a)_{a\in Q}\\
K(x_a,x_{a^\ast})_{a\in
  Q}&=(-ix^\dagger_{a^\ast},ix^\dagger_a)_{a\in Q}
\end{align*}
It is clear that with respect to this quaternionic structure the
metric \eqref{ref-2.4-16} is hyper-K\"ahler. Let $\HH^0$ be the kernel of
the reduced trace map on $\HH$. If $\beta\in\HH^0$ then there is an
associated real symplectic form on $V$ defined by
$\omega_\beta(v,w)=(v,\beta w)$. 

Let us write $\mathfrak{gl}=\Lie(\Gl(\alpha))$ and
$\mathfrak{u}=\Lie(U(\alpha))$ where $U(\alpha)$ is the maximal compact
subgroup of $\Gl(\alpha)$ given by the product of unitary groups
$\prod_{i\in I}
U(\alpha_i)$. The hyper-K\"ahler structure on $V$ is clearly
$U(\alpha)$-invariant and it is a standard fact  that the symplectic form
$\omega_\beta$ has an associated  moment map $\mu_\beta:V\r \mathfrak{u}^\ast$
 given by
$\mu_\beta(v)(u)=-\frac{1}{2}\omega_\beta(v, uv)$ for $x\in V$,
$u\in\mathfrak{u}$. Below we will write $\mu_\RR$ for $\mu_I$.

The three moment maps $\mu_I$, $\mu_J$, $\mu_K$ may be combined into a
so-called hyper-K\"ahler moment map 
\begin{equation}
\label{ref-2.5-17}
\mu:V\r \HH^0\otimes_\RR \mathfrak{u}^\ast: x\mapsto I\otimes \mu_I(x)+J\otimes
\mu_J(x)+K\otimes \mu_K(x) 
\end{equation}
{}From the explicit description of $\mu_\beta$ we deduce for
$h\in\HH$:
\begin{equation}
\label{ref-2.6-18}
\mu_\beta(hx)=\mu_{\bar{h}\beta h}(x)
\end{equation}
where $\bar{h}$ is the conjugate of $h$ in $\HH$.  {}From
\eqref{ref-2.6-18} we deduce that \eqref{ref-2.5-17} is
$\HH^\ast$-invariant if we let $\HH^\ast$ act on $\HH^0$ by $h\cdot
\beta=h\beta \bar{h}$. 

For this action $\HH^0-\{0\}$ is a homogeneous space and hence if we
choose $\beta\in \HH^0-\{0\}$ and a contractible subset $S\subset
\HH^0-\{0\}$ containing $\beta$ then there is a continuous map
$\theta_{\beta,S}:S\r \HH^\ast$ which is a section (above $S$) for the
map $h\mapsto h\cdot \beta$. 

Choose a $U(\alpha)$-invariant $\lambda\in \mathfrak{u}^\ast$ and let $V'=\mu^{-1}(S\times
\lambda)$, $V''=\mu^{-1}(\beta\times \lambda)$. Then $V''\times
S\r V':(x,s)\mapsto \theta_{\beta,S}(s)x$ defines a trivialization of
$\mu\mid V''$. Thus we have proved that above $S\times \lambda$, $\mu$
is a trivial bundle. Moreover this trivialization is clearly
$U(\alpha)$-equivariant.

Put $\omega_\CC=\omega_J+i\omega_K$. This is a complex
$\Gl(\alpha)$-invariant symplectic form on $V$ and it is easy to see that the
associated moment map $V\r \mathfrak{gl}^\ast$ is given by
$\mu_\CC(x)=\mu_J(x)+i\mu_K(x)$ where we have extended
$\mu_J(x),\mu_K(x)$ to linear maps $\mathfrak{gl}\r \CC$. A straightforward
computation shows that 
\begin{align*}
\mu_\RR(x)&=\frac{i}{2}\sum_a [x_a,x^\dagger_a]\\
\mu_\CC(x)&=\sum_{a\in Q} [x_a,x_{a^\ast}]
\end{align*}
where we have identified $\mathfrak{u}$, $\mathfrak{gl}$ with their duals via  the
trace form  $(g,h)=-\Tr(gh)$ (the minus sign makes the form positive
definite on $\mathfrak{u}$).

{}From the
description $\mu_\CC=\mu_J+i\mu_K$ we obtain:
\[
\mu_\CC^{-1}(a)=\mu_J^{-1}\left(\frac{a-a^\dagger}{2}\right) \cap\mu_K^{-1}\left(\frac{a+a^\dagger}{2i}\right)  
\]
which yields
\begin{align*}
\mu^{-1}_\CC(\CC\lambda)\cap \mu^{-1}_\RR(i\lambda)&\cong \mu^{-1}((I+\RR
J +\RR K)\times i\lambda)\\
\mu^{-1}_\CC(0)\cap \mu^{-1}_\RR(i\lambda)&\cong \mu^{-1}(I\times i\lambda)\\
\end{align*}
{}From the fact that $I+\RR
J +\RR K$ is contractible  we deduce as explained above that $\mu$ is
trivial above $(I+\RR
J +\RR K)\times i\lambda$. Since on the inverse image of $(I+\RR
J +\RR K)\times i\lambda$, $\mu$ and $\mu_\CC$ are basically the same
we deduce that $\mu_\CC:\mu_\CC^{-1}(\CC\lambda)\cap
\mu_\RR^{-1}(i\lambda)\r \CC\lambda$ is a trivial family in a way
that is compatible with the $U(\alpha)$-action.

We now use this to construct the following commutative diagram of
continuous maps:
\[
\begin{CD}
X_s\times L @>\pr>> L\\
 @ArAA @|\\
\bigl(\mu^{-1}_\CC(0) \cap
\mu^{-1}_\RR(i\lambda)\bigr)/U(\alpha)\times L
@>\pr>>L\\
@Vp  VV
@|\\
 \bigl(\mu^{-1}_\CC(L)
\cap \mu^{-1}_\RR(i\lambda)\bigr)/U(\alpha)
@>\bar{\mu}_\CC>> L\\
@Vr' VV @|
\\
\Xi
@>f>> L
\end{CD}
\]
Here $p$ is obtained from the trivialization of $\mu_\CC$ we have 
constructed above (recall that $L=\CC\lambda$) and $r,r'$ are obtained from the inclusion
$\mu^{-1}_\RR(i\lambda)\subset \Rep(\bar{Q},\alpha)^\lambda$
\cite[Prop. 6.5]{King}.

To prove the lemma it is now sufficient to show that the vertical maps
on the left are homeomorphisms.
This is true by construction for $p$. We claim that it is also true
for $r,r'$. It suffices to consider $r'$ since $r$ is obtained from
$r'$ by restricting to a fiber.

By \cite[Prop. 6.5]{King} $r'$ is a bijection. Hence it suffices to
show that $r'$ is proper. Clearly $r'$ is the restriction to 
$ \bigl(\mu^{-1}_\CC(L)
\cap \mu^{-1}_\RR(i\lambda)\bigr)/U(\alpha)$ of the first map in the
following diagram
\[
\mu^{-1}_\RR(i\lambda)/U(\alpha)\r \Rep(\bar{Q},\alpha)^\lambda\quot G(\alpha)
\r
\Rep(\bar{Q},\alpha)\quot G(\alpha)
\]
By Theorem \ref{ref-2.3.4-19} below  the composition of these two maps
is proper. It follows that 
the first map is also proper. This finishes the proof.
\end{proof}
We have used the following result.
\begin{theorems} 
\label{ref-2.3.4-19}
\cite[Theorem 1.1]{Neeman2} Let the notations be as above.
The canonical map
\[
\psi:V\r V\quot G \times \mathfrak{u}:v\mapsto (\bar{v}, \mu_\RR(v))
\]
is proper. 
\end{theorems}

\subsection{End of proof}
Let $k=\overline{\FF}_p$. We choose $\lambda$ generic with respect to
$\alpha$. Now recall that Kac has shown \cite{kac2} that $a_\alpha(q)$
is a polynomial. We first 
show that $X_k$ is pure. By Corollary \ref{ref-2.3.2-14} we may as
well show that $X_{s,k}$ is pure. 
\emph{Since we will now  work exclusively over $k$ we  drop
the corresponding subscript.}

Define $X^0_s=\Rep(\Pi^0,\alpha)\quot G(\alpha)$. Then the
canonical map $u:X_s\r X^0_s$ is projective \cite{King}. Let
$\GG_m$ act on $\Rep(\bar{Q},\alpha)$ in such a way that $\eta\in
\GG_m$ multiplies all arrows by $\eta$. This action induces
$\GG_m$-actions on $X_s$ and $X^0_s$ and the map $u$
commutes with these actions. 

Now clearly $X^0_s=\Spec R$  with
$R=\Oscr(\Rep(\Pi^0,\alpha))^{G(\alpha)}$.
The ring $R$ is  graded 
via the $\GG_m$-action we have defined in the
previous paragraph and it is easy to see that the grading is of the
form $R=k+R_1+R_2+\cdots$ with $R_i$ finite dimensional.

Thus it follows that $(X^{0}_s)^{\GG_m}_0$ consists
of a single point $o$
defined by the  graded maximal ideal of $R$. It also follows that
$(X_s)^{\GG_m}\subset u^{-1}(o)$. Since  $u$ is
projective it follows that $(X_s)^{\GG_m}$ is also projective. Hence by
Proposition \ref{ref-A.2-26} $X_s$ is pure. 

By combining Proposition \ref{ref-2.2.1-12}, Lemma \ref{ref-A.1-24} with
Corollary \ref{ref-2.3.2-14} it follows
\[
a_\alpha(q)=\sum_{i\ge 0}
\dim H^{2d+2i}_c(X_{s,k},\QQ_l)q^i
\]
with $d=1-(\alpha,\alpha)$ and $k=\overline{\FF}_p$ for $p\gg
0$.  Since this is true for large characteristic we obtain
\begin{equation}
\label{ref-2.7-20}
a_\alpha(q)=\sum_{i\ge 0}
\dim H^{2d+2i}_c(X_{s,\CC},\CC)q^i
\end{equation}
Furthermore if $X_{s,\CC}$ is non-empty then we compute 
\[
\dim
X_{s,\CC}=\dim \Rep(\Pi^0,\alpha)^\lambda-\dim G(\alpha)=
\dim \Rep(\bar{Q},\alpha)^\lambda-2\dim G(\alpha)=2d
\]
Thus the sum in \eqref{ref-2.7-20} runs from $i=0$ to $i=d$. Applying
Poincar\'e duality we obtain \eqref{ref-1.1-2}.

Finally we prove the equivalence of Conjecture B and Proposition \ref{ref-1.5-3}.
\emph{In the rest of this section our base field will be $\CC$}. 
 Our starting point is
the following commutative diagram 
\begin{equation}
\label{ref-2.8-21}
\begin{CD}
\Rep(\Pi^0,\alpha)^\lambda @>\text{open} >>
\Rep(\Pi^0,\alpha)\\
@VVV @VVV\\
\Rep(\Pi^0,\alpha)^\lambda/G(\alpha) @>>u>
\Rep(\Pi^0,\alpha)\quot G(\alpha)
\end{CD}
\end{equation}
where all the maps are the obvious ones. 

By \eqref{ref-1.1-2} we have
$a_\alpha(0)=\dim H^{2d}(X_s,\CC)$. 
With a similar argument as the one used 
in \cite[Prop. 4.3.1]{slodowy}  one shows that $X_s$ is homotopy 
equivalent to $u^{-1}(0)$. Thus $H^{2d}(X_s,\CC)=H^{2d}(u^{-1}(0),\CC)$.

Let $(-)^{\text{nil}}$ denote the nilpotent
representations in $\Rep(\Pi^0,\alpha)$ and $\Rep(\Pi^0,\alpha)^\lambda$.
Since the leftmost map in
\eqref{ref-2.8-21} is a principal $G(\alpha)$-bundle and the top map is
an open immersion we find that if $X_s\neq \emptyset$ then $\dim
u^{-1}(0)=\dim \Rep(\Pi^0,\alpha)^{\text{nil}}-\dim G(\alpha)$. Since
$\Rep(\Pi^0,\alpha)^{\text{nil}}$  \cite[12.9]{Lusztig2}  
is a Langrangian subvariety of $\Rep(\bar{Q},\alpha)$. Thus
$u^{-1}(0)$ is equidimensional and furthermore 
$\dim
u^{-1}(0)=(1/2)\dim \Rep(\bar{Q},\alpha)-\dim G(\alpha)=d$. Hence (even
if $X_s=\emptyset$), $\dim H^{2d}(X_s,\CC)$ is equal to the number of
irreducible components of $u^{-1}(0)$. Using again that the leftmost
map is a principal $G(\alpha)$-bundle this is equal to the number of
irreducible components of
$\Rep(\Pi^0,\alpha)^{\lambda,\text{nil}}$. This finishes the proof.
\section{Proof of Conjecture B for indivisible roots}
\label{newsection}
\emph{In this section  our ground field is $\CC$.} 

At the end of the previous section it was shown that Proposition
  \ref{ref-1.5-3} and Conjecture B are equivalent. So we only prove
  Proposition  \ref{ref-1.5-3}. The idea for the proof of Proposition
  \ref{ref-1.5-3} 
  came partially from a talk by Reineke \cite{Reineke}.

In the previous section we have used the notion of $\lambda$-stability
introduced by King \cite{King} which is derived from geometric
invariant theory. A technical inconvenience of this notion is that if
we work in $\Rep(\bar{Q},\alpha)$ then $\lambda\cdot \alpha$ must be
zero. Hence we cannot use the same $\lambda$ for all
$\alpha$. Following Reineke \cite{Reineke} we use therefore an alternative
notion of stability we will call slope stability.

We fix an element $\Theta\in\ZZ^I$ and we define the
corresponding ``slope function'' $s(\alpha)=(\Theta\cdot\alpha)/ \dim
\alpha$ where $\dim\alpha=\sum\alpha_i$.  If $V$ is a finite
dimensional representation of $\bar{Q}$ then we put
$s(V)=s(\underline{\dim} V)$.  If $X\subset \Rep(\bar{Q},\alpha)$ is
irreducible then we write $s(X)$ for the slope of a generic point of
$X$.

 A representation $V$ of
$\bar{Q}$ is ($\Theta$-slope) stable (resp. semistable) if for all proper
subrepresentations $W$ of $V$ we have $s(W)<s(V)$ (resp. $s(W)\le
s(V)$). 
It is easy to see that for a fixed dimension vector $\alpha$, King (semi)stability and
slope (semi)stability are equivalent for suitable $\lambda$ and
$\Theta$. Below the notion of (semi)stability will refer to
$\Theta$-slope (semi)stability for an arbitrary but fixed $\Theta$.

The following lemma is standard.
\begin{lemma} \label{classical} Assume that $V,W$ are semistable representations such
  that $s(V)>s(W)$. Then $\Hom(V,W)=0$.
\end{lemma}
The following result is proved in \cite{HN,Reineke}.
\begin{theorem}
\label{maintheorem} Let $V$ be a representation of $\bar{Q}$. Then there
  exists a unique filtration 
\[
0=V_0\subsetneq V_1\subsetneq\cdots \subsetneq V_{n-1}
\subsetneq V_n=V
\]
such that all $V_{i+1}/V_i$ are semistable and such that
$s(V_{i+1}/V_i)$ is a strictly decreasing function of $i$. 
\end{theorem}
The filtration introduced in the last theorem is called the
Harder-Narasimhan filtration. Let us write
\[
t(V)=(\udim (V_1/V_0),\ldots, \udim (V_n/V_{n-1}))
\]
We call $t(V)$ the HN-type of $V$.

If $X$ is a variety then we write $\Irr X$ for the set of irreducible
components of $X$. If $\alpha\in\NN^I$ then we write $\Lambda_\alpha$ for
$\Rep(\Pi^0,\alpha)^{\text{nil}}$. According to \cite{Lusztig2} this is a
Lagrangian subvariety  of $\Rep(\bar{Q},\alpha)$ and furthermore
$\Irr \Lambda_\alpha$ indexes a basis for $U(\frak{g}^+)_\alpha$\cite{SK,Lusztig5}.

If $X\subset \Irr\Rep(\bar{Q},\alpha)$ we say that $X$ is semistable if it
contains a semistable representation.  We write $s(X)$
for $s(V)$ with $V\in X$ generic.

Let $S_\alpha$ be the set of tuples $Z^\ast=(Z_1,\ldots,Z_{n})$  with $Z_i$ semistable elements of certain $\Irr
\Lambda_{\alpha_i}$ such that 
$\alpha=\sum \alpha_i$  and such that
$s(Z_i)$  is strictly decreasing.

For  $Z^\ast\in S_\alpha$ we define $m'(Z^\ast)$ as the
set of all $V\in \Lambda_\alpha$ such that if $(V_i)_i$ is the
HN-filtration on $V$ then $V_{i}/V_{i-1}\in Z_i$. 

The following is our main theorem.
\begin{theorem}  
\begin{enumerate}
\item
If $Z^\ast\in S_\alpha$ then $m'(Z^\ast)$ has a dense intersection with
unique $Z \in \Irr\Lambda_\alpha$. Put $m(Z^\ast)=Z$.
\item
The map $m$ defines a bijection between  $S_\alpha$  and $\Irr \Lambda_\alpha$.
\end{enumerate}
\end{theorem}
\begin{proof}
By the existence and uniqueness of the HN-filtration 2. follows from
1. Hence we only have to prove 1.

Let us call a subset $Z$ of $\Lambda_\alpha$ \emph{good} if it has the following
properties.
\begin{enumerate}
\item The elements of $Z$ have constant HN-type.
\item $Z$ is constructible. 
\item $Z$ has a dense intersection with a unique irreducible
  component of $\Lambda_\alpha$.
\end{enumerate}

By induction it is clearly sufficient to prove the following claim:

\smallskip

\noindent \textbf{Claim.} Let $Z_1\in \Irr \Lambda_\beta$ be
semistable
and let  $Z_2\subset \Lambda_\gamma$ be good.
Assume that $s(Z_1)_1>t(V)_1$ for $V\in Z_2$ arbitrary.
Define
  $Z\subset\Lambda_{\beta+\delta}$ as the set of all $V\in
   \Lambda_{\beta+\delta}$ which contain a semistable subrepresentation $U\subset
  V$ such that $U\in Z_1$, $V/U\in Z_2$. Then $Z$ is good.

\smallskip
The only non-obvious property to prove is that $Z$ has a dense
intersection
with a unique irreducible component of $\Lambda_{\beta+\gamma}$. So
this is what we do below.

Let $Z^\circ_1$ be the semi-stable locus of $Z_1$ and let $E$ be the
set of 5-tuples $(U,V,W,u,w)$ with $U\in Z^\circ_1$, $V\in
\Lambda_{\beta+\gamma}$, $W\in Z_2$, $u\in \Hom(U,V)$, $w\in\Hom(V,W)$
such that
\[
0\r U\xrightarrow{u} V\xrightarrow{w} W\r 0
\]
is exact. It is easy to see  that $E$ is a constructible
subset of $\Rep(Q,\beta)\times \Rep(Q,\beta+\gamma)\times
\Rep(Q,\gamma)\times \Hom_k(U,V)\times \Hom_k(V,W)$. 

Due to the uniqueness of the HN-filtration the non-empty fibres of the
projection map $p:E\r \Lambda_{\beta+\gamma}:(U,V,W,u,w)\mapsto V$ are
isomorphic to $\Gl(\alpha)\times \Gl(\beta)$ and hence they have
dimension $\alpha\cdot\alpha+\beta\cdot\beta$.

There is another projection map $q:E\r Z^\circ_1\times
Z_2:(U,V,W,u,w)\mapsto (U,W)$. According to
\cite[Lemma 5.1]{CB2} its  fibers have dimension
\[
(\beta+\gamma)\cdot (\beta+\gamma)+\dim \Ext^1(W,U)-\dim \Hom(W,U)
\]
and the proof also shows that these fibers are irreducible and locally
closed. 

According to \cite[Lemma 1]{CB1} we also have
\[
\dim \Hom(U,W)-\Ext^1(W,U)+\dim \Hom(W,U)=2(\beta,\gamma)
\]
and furthermore  according to lemma \ref{classical} we have $\Hom(U,W)=0$. 
Substituting we find that the fibers of $q$ have dimension:
\[
(\beta+\gamma)\cdot (\beta+\gamma)-2(\beta,\gamma)
\]
According to lemma \ref{elementary} below we find that 
$E$  contains a dense irreducible locally closed subset $E'$  such that $\dim
(E-E')<\dim E$. Furthermore the dimension of $E$ is:
\begin{equation}
\label{firstdimension}
\dim \Lambda_\beta+\dim \Lambda_\gamma+(\beta+\gamma)\cdot (\beta+\gamma)-2(\beta,\gamma)
\end{equation}
Now we have for $\alpha\in\ZZ^I$:
\[
\dim \Lambda_\alpha=\frac{1}{2}\dim\Rep(\bar{Q},\alpha)=\alpha\cdot
\alpha-(\alpha,\alpha)
\]

A trite computation shows that $Z=p(E)$ has dimension
\[
(\beta+\gamma)\cdot
(\beta+\gamma)-(\beta+\gamma,\beta+\gamma)=\dim
\Lambda_{\beta+\gamma}
\]
 and 
$p(E-E')$ has smaller dimension.
Hence $\dim p(E')=\dim \Lambda_{\beta+\gamma}$.
Since $E'$ is irreducible it follows that $p(E')$ is dense in some
irreducible component $Z$ of $\Lambda_{\beta+\gamma}$. This finishes
the proof.
\end{proof}
If $X$ is an algebraic variety and $S\subset X$ is a constructible set
then let us say that $S$ is weakly irreducible if $S$ contains a
dense subset $S'$ which is irreducible locally closed in $X$ and has
the property that
$\dim(S-S')<\dim S$.
\begin{lemma}
\label{elementary}
  Let $q:X\r Y$ be a morphism between (reduced) algebraic varieties.
  Let $S\subset X$, $T\subset Y$ be constructible subsets with
  $T=q(S)$ such that the fibers of $q:S\r T$ are locally closed
  in $X$, irreducible and of constant dimension. If $T$ is weakly
  irreducible then so is $S$.
\end{lemma}
\begin{proof}
Left to the reader.
\end{proof}
For $\alpha\in \NN^I$ let us put $n_\alpha$ for the number of
components of $\Lambda_\alpha$ and $m_\alpha$ for the number of
semistable components.  By \cite{Lusztig5} we have $n_\alpha=\dim
U(\frak{g}^+)_\alpha$. Theorem \ref{maintheorem} yields the formula
\[
n_\alpha=\sum_{\begin{smallmatrix}\alpha_1,\ldots,\alpha_n
\\
s(\alpha_1)>\cdots >s(\alpha_n)
\\ \sum \alpha_i=\alpha, 
\end{smallmatrix}
}
\prod_i m_{\alpha_i}
\] 
and this formula allows us to determine the $m_\alpha$ recursively
from the $n_\alpha$. 

Put $r_\alpha=\dim \frak{g}_\alpha$.
It turns out  that we can give an explicit expression for the $m_\alpha$ in
terms of $r_\alpha$.  Put an arbitrary
total ordering on $\NN^i$ with the property $s(\beta)>s(\gamma)\Rightarrow
\beta> \gamma$ and $\beta>\gamma\Rightarrow s(\beta)\ge
s(\gamma)$. 
\begin{lemma} The following formula holds.
\begin{equation}
\label{maformula}
m_\alpha=\sum_{\begin{smallmatrix}(u_1,\beta_1),\ldots,(u_n,\beta_n)
\\ \beta_1>\ldots>\beta_n
\\ s(\beta_1)=\cdots=s(\beta_n)=s(\alpha)
\\ \sum u_i\beta_i=\alpha, 
\end{smallmatrix}
}
\prod_i  \choose{r_{\beta_i}+u_i-1}{u_i}
\end{equation}
\end{lemma}
\begin{proof}
By the PBW-theorem we
have
\[
n_\alpha=
\sum_{\begin{smallmatrix}(u_1,\beta_1),\ldots,(u_n,\beta_n)
\\ \beta_1>\ldots>\beta_n
\\ \sum u_i\beta_i=\alpha, 
\end{smallmatrix}
}
\prod_i  \choose{r_{\beta_i}+u_i-1}{u_i}
\]
In this formula we may collect the $\beta_i$'s with equal
slope. Let $m'_\alpha$ be given by the righthand side of
\eqref{maformula}.
Then we have
\[
n_\alpha=\sum_{\begin{smallmatrix}\alpha_1,\ldots,\alpha_n
\\
s(\alpha_1)>\cdots >s(\alpha_n)
\\ \sum \alpha_i=\alpha, 
\end{smallmatrix}
}
\prod_i m'_{\alpha_i}
\] 
and by induction it follows $m'_{\alpha_i}=m_\alpha$. This finishes the
proof of \eqref{maformula}.
\end{proof}
\begin{proof}[Proof of Proposition \ref{ref-1.5-3}]
  Recall that $\lambda\in \ZZ^I$ is such that $\lambda\cdot \alpha=0$ and
  $\lambda\cdot \beta\neq 0$ for all $0<\beta<\alpha$.

 Now it is clear that
King semistability for $\lambda$ is equivalent to 
slope semistability for $\Theta=-\lambda$. Hence for this particular
$\Theta$ we need to show that $m_\alpha=r_\alpha$. This follows
immediately from \eqref{maformula}.
\end{proof}

\appendix
\section{Purity}
\label{ref-A-23}
For the benefit of the reader we recollect some basics. As usual $q$ is a power of a prime number
$p$ and $l\neq p$
is another prime number. We put $k=\bar{\FF}_p$.

Assume that $Z/k$ is a variety defined over
$\FF_q$, i.e. there is some $Z_0/\FF_q$ such that
$Z=(Z_0)_{k}$. Let ${{F}}:Z\r Z$ be the corresponding Frobenius
morphism. The key method for counting rational points on $Z_0$ is given by the
trace formula \cite[Thm
  3.2]{Rapport}
\[
|Z_0(\FF_{q^r})|=\sum^{2\dim Z}_{i=0} (-1)^i \Tr({{F}}^r; H^i_c(Z,\QQ_l))
\]
For this formula to be effective one needs information on the
eigenvalues of ${{F}}$. Let us say
that $Z$ is \emph{(cohomologically) pure} if the eigenvalues of ${{F}}$ acting on
$H^i_c(Z,\QQ_l)$  have absolute value $q^{i/2}$. This
definition only depends on $Z$ and not on the particular choice of
$\FF_q$ and $Z_0$.  The Weil
conjectures
\cite{De} imply that if $Z$ is smooth proper over $k$ then $Z$ is
pure. 

We have used the notion of purity in the following context: 
\begin{lemma} 
\label{ref-A.1-24}
Assume that $Z$ is pure and that there is a polynomial
  $p(t)\in\ZZ[t]$ such that $|Z_0(\FF_{q^r})|=p(q^r)$. Then
  $p(q^r)=\sum_i \dim H^{2i}_c(Z,\QQ_l) q^{ri}$ and in particular
  $p(t)\in\NN[t]$. 
\end{lemma}
\begin{proof}
It is clearly sufficient to show that the action of ${{F}}$ on $H^{2i}_c(Z,\QQ_l)$
has a unique eigenvalue $q^{i}$ and that in addition $H^{2i+1}_c(Z,\QQ_l)=0$.

Write $p(t)=\sum_i b_{2i} t^i$ and $b_j=0$ for $j$ odd. Since $Z$ is
pure the eigenvalues of ${{F}}$ acting on $H^i(Z,\QQ_l)$ are given by
$\epsilon_{ij} q^{i/2}$ where $j=1\ldots \beta_i$ and
$|\epsilon_{ij}|=1$. {}From the hypotheses and the trace formula we obtain
\begin{equation}
\label{ref-A.1-25}
\sum_{i=0}^{2d} (-1)^i b_i q^{ri/2}=\sum_{i=0}^{2d} (-1)^i \sum_{j=1}^{\beta_i}
\epsilon_{ij}^r q^{ri/2}
\end{equation}
where $d=\dim Z$.
Dividing by $q^{rd}$ we find 
\[
b_{2d}=\lim_{r\r \infty} \sum_{j=1}^{\beta_{2d}}
\epsilon_{2d,j}^r
\]
Using a Van der Monde type argument we see that the limit on the 
righthand side only exists if  $\epsilon_{2d,j}=1$ for all
$j$. Subtracting the leading term in $q$ from \eqref{ref-A.1-25} and
repeating the same argument we
ultimately find that  $\epsilon_{ij}=1$ for all $i,j$. Since $b_i=0$
for odd $i$ we find that $\beta_i=0$ for odd $i$. This finishes the
proof. 
\end{proof}
In this paper we  use the following purity criterion:
\begin{proposition} 
\label{ref-A.2-26}
Assume that $Z$ is smooth quasi-projective and that there is an action $\lambda:\GG_m\times
Z\r Z$ such that for every $x\in Z$ the limit
$\lim_{t\r 0}\lambda(t,x)$ exists. Assume in addition that
$Z^{\mathbb{G}_m}$ is projective. Then $Z$ is pure. 
\end{proposition}
\begin{proof}
Let $Z^{G_m}=\bigcup_\alpha L_\alpha$ be the decomposition
into connected components and for each $\alpha$ define
\[
W_\alpha=\{x\in Z\mid \lim_{t\r 0} \lambda(t,x)\in L_\alpha\}
\]
According to \cite[Thm 4.1, proof of Thm 4.2] {Bia} the $L_\alpha$,
$W_\alpha$ are smooth and the $W_\alpha$ are locally closed in
$Z$. Furthermore the limit map $f_\alpha:W_\alpha\r L_\alpha$ is a Zariski
locally trivial affine fibration.
 Furthermore in \cite{Bia1} it is shown that there is a filtration
 $\emptyset=Z_0\subset Z_1\subset \cdots \subset Z_n=Z$ of $Z$ by closed
 subsets such that for every $i$, $Z_{i+1}-Z_i$ is one of the
 $W_\alpha$ (this depends on $Z$ being quasi-projective).

Looking at Zariski open sets we  find
\[
R^if_{\alpha\ast} \QQ_l=
\begin{cases} 
\QQ_l&\text{if $i=0$}\\
0&\text{otherwise}
\end{cases}
\]
Thus
\begin{equation}
\label{ref-A.2-27}
H^i(W_\alpha,\QQ_l)=H^i(L_\alpha,\QQ_l)
\end{equation}
By the Weil conjectures $L_\alpha$ is pure. 
Since $L_\alpha$ and $W_\alpha$ are smooth,  \eqref{ref-A.2-27} and  lemma \ref{ref-A.3-28} below
imply that $W_\alpha$ is smooth as well. Applying lemma \ref{ref-A.4-29}
finishes the proof.
\end{proof}
We have used the following lemmas
\begin{lemma}
\label{ref-A.3-28}
If $Z$ is smooth then $Z$ is pure if and only if
the  eigenvalues of ${{F}}$  acting on
$H^i(Z,\QQ_l)$ have absolute values $q^{i/2}$.
\end{lemma}
\begin{proof}
This follows by Poincar\'e duality.
\end{proof}
\begin{lemma}
\label{ref-A.4-29}
 Assume that we have a decomposition $Z=Y\coprod U$
  where
$Y$ is closed and $Y$, $U$ are pure. Then $Z$ is also pure and
  in addition we have short exact sequences
\begin{equation}
\label{ref-A.3-30}
0\r H^i_c(Y,\QQ_l)\r H^i_c(Z,\QQ_l)\r H^i_c(U,\QQ_l)\r 0
\end{equation}
\end{lemma}
\begin{proof}
This follows from  the fact that in the long exact sequence
\[
\r H^{i-1}_c(U,\QQ_l)\r H^i_c(Y,\QQ_l) \r H^i_c(Z,\QQ_l)\r
H^i_c(U,\QQ_l)\r
H^{i+1}_c(Y,\QQ_l)\r
\]
the connection maps must be zero by purity.
\end{proof}

\section{Invariants over $\ZZ$}
\label{ref-B-31}
In this paper we have used lifting to characteristic zero. To do this
rigorously we need that taking invariants commutes with base change
over a Zariski open part of the base. This is of course well known 
but we have not found an
explicit reference.
For simplicity we will only consider the case where the base is
$\ZZ_f$. Replacing $\Spec \ZZ_f$ by a Zariski open subset amounts to
``increasing'' $f$ in the following sense:
\begin{convention} If $f\in\ZZ$ then \emph{increasing} $f$ means
  making $f$ larger for the partial order given by divisibility. 
\end{convention}

Let $G$ be reductive
group defined over  $\ZZ_f$ \cite{seshadri}. All $G$-actions below are
rational. That is: they are obtained from a coaction of $\Oscr(G)$.

First recall Seshadri's generalization of Geometric Invariant Theory
to an arbitrary base ring.
\begin{theorem}\cite[\S II]{seshadri}
Let $R$ be finitely generated $\ZZ_f$ algebra  and let $M$ be a
finitely generated $R$-module. Assume that $G$ acts rationally 
on $R$ and $M$. Then $R^G$ is a finitely
generated $\ZZ_f$-algebra and $M^G$ is a finitely generated $R^G$-module. In addition 
if $X=\Spec R$ and $X\quot G=\Spec R^G$ then $X\quot G$ has the usual behavior in the sense that if
$\Spec k\r \Spec\ZZ_f$ is a geometric point
then  the points in $(X\quot G)(k)$ correspond to the closed orbits in
$X(k)$. 
\end{theorem}
It follows in particular that  $\Spec (R\otimes
k)^G\r \Spec (R^G\otimes k)$ is set-theoretically a bijection. We
want it to be an isomorphism.
The result we need is the following:
\begin{theorem}
\label{ref-B.3-32}Let $R$ be finitely generated $\ZZ_f$ algebra and let $M$ be a
finitely generated $R$-module. Assume that $G$ acts rationally 
on $R$ and $M$. Then there exists a Zariski open subset $U$ of $\Spec \ZZ_f$
  such that for every geometric point $\Spec k\r U$ we have that the
  canonical map $M^G\otimes k\r (M\otimes k)^G$ is an isomorphism and
  in addition $H^i(G,M\otimes k)=0$ for $i>0$.
\end{theorem}
We will informally say that the formation of $M^G$ is compatible with
base change for $f$ large enough.
\begin{proof} 
Recall that if $H$ is a reductive algebraic group over an
algebraically closed field $k$  then an $H$-representation of
countable dimension is said to have a good
filtration if it has an ascending filtration by co-Weyl-modules $Y(\lambda)$, or
equivalently if $H^i(H,Y(\lambda)\otimes U)=0$ for all $i>0$ and all $\lambda$ \cite{Donkin2}. In particular $(-)^H$ is exact on
representations with a good filtration and the category of
representations with good filtrations is stable under taking
cokernels of  surjective 
maps and extensions.
It is a deep theorem \cite{Donkin2,Mathieu} that the category
of representations with a good filtration is stable under tensor
product.

Put $A=\ZZ_f$.
If $V$ is a $G$-module free of finite rank over $A$ and if $V\otimes_A
k$ ($k$ as in the statement of the theorem) has a good filtration then it follows from exactness of
$(-)^G$ that $\dim (V\otimes_A k)^G$ is the number of
$Y(0)$'s in a good filtration of $V\otimes_A k$. This can be
computed in terms of characters so we conclude 
\begin{equation}
\label{ref-B.1-33}
\dim (V\otimes_A k)^G=\dim(V\otimes_A \bar{\QQ})^G=\rk V^G=\dim
(V^G\otimes_A k)
\end{equation}
By the universal coefficient theorem the canonical map 
\begin{equation}
\label{ref-B.2-34}
V^G\otimes_A k\r (V\otimes_A k)^G
\end{equation}
is a monomorphism  and hence by \eqref{ref-B.1-33} it an isomorphism. 

If $V$ is not necessarily of finite rank but has a filtration
\hbox{$0=V_0\subset V_1\subset V_2\subset\cdots$} such that each
$V_{i+1}/V_i$ is free of finite rank and $(V_{i+1}/V_i)\otimes_A k$ has a good
filtration then it is easy to see that \eqref{ref-B.2-34} is still an
isomorphism.

  Since the action of $G$ is locally finite there exist a
finitely generated $G$ module $W$ such that $R$ is a quotient
of $SW$. By increasing $f$ we may assume that $W$ is free. If the characteristic of $k$ is large with respect to $\lambda$ (in
a suitable sense) then $Y(\lambda)$ is  simple \cite[Ch. 6]{Jantzen}.
It follows that if $\kar k$ is large then the finite
dimensional $G$-representation $\Lambda(W\otimes_A k)$ has a good
filtration. It then follows from \cite[\S 4.3]{JA} that
$SW\otimes_A k=S(W\otimes_{A} k)$ has a good filtration as well. {}From
the proof it follows that this good filtration is compatible with the grading.

Now we filter $SW$ by degree and we put
the induced filtration on $R$. We choose a compatible filtration
on $M$ such that $\gr M$ is a finitely generated $\gr
R$-module (confusingly such a filtration is also called a good
filtration!) \cite{NVO}.  Since $\gr R$ and $\gr M$ are finite over the
noetherian ring $SV$ their $\ZZ$-torsion is supported
  on a finite set of primes. Hence by increasing $f$ we may and we
  will assume
  that $\gr R$ and $\gr M$ are torsion free.

Since $SW$ has finite global dimension it is easy to see that (at the
cost of possibly increasing $f$) we may construct a graded resolution
of $\gr M$ whose terms are of the form $U_i\otimes_A SW$ with $U_i$ a  free
$G$-representation of finite
rank. Increasing $f$ again if necessary we may assume that all
$U_i\otimes_A k$ have a good filtration. Thus it follows that $(\gr M)\otimes_A k$ will
also have a good filtration compatible with the grading for all
$k$. Thus $M\otimes_A k$ has vanishing cohomology. The rest of the
theorem follows from
the fact that \eqref{ref-B.2-34} is an isomorphism with $V=M$.
\end{proof}
{}From Theorem \ref{ref-B.3-32} one easily deduces that all standard
constructions are compatible with base change if we take $f$ large
enough. We give an  example whose proof we leave to the reader.
\begin{lemma} Let $X$ be of finite type over $\ZZ_f$ and assume that
  $G$ acts rationally on $X$. Let $L$ be a $G$-equivariant line bundle on $X$.
  Let $X^{ss}$ be the $L$-semistable points on $X$ \cite[\S
  II]{seshadri}. Then the formation of $X^{ss}$ and $X^{ss}\quot G$
  is compatible with base change for $f$ large enough.
\end{lemma}

\section*{Appendix by Hiraku Nakajima}

The following simple proof avoids the arguments in Section~\ref{cohomsection},
showing directly that if $\lambda$ is generic for $\alpha$, then
$\operatorname {Rep}(\Pi^\lambda,\alpha)^\lambda/\!\!/ G(\alpha)$
and $\operatorname {Rep}(\Pi^0,\alpha)^\lambda/\!\!/ G(\alpha)$ have the same
number of points over sufficiently large finite fields $\mathbb{F}_{q}$. 

Let $k = \overline{\mathbb F_p}$, the algebraic closure of a finite
field.

Suppose that $\pi\colon \mathcal X \to \mathbb A^1$ is a smooth family
of nonsingular quasi-projective varieties over the line 
$\mathbb A^1 = k$ with the following properties:
\begin{enumerate}
\item there exists a $\mathbb G_m$-action on $\mathcal X$ such that
$\pi$ is equivariant with respect to a $\mathbb G_m$-action on
$\mathbb A^1$ of positive weight,
\item for every $x\in\mathcal X$, the limit $\lim_{t\to 0} t\cdot x$
exists.
\end{enumerate}

Let $X_\lambda = \pi^{-1}(\lambda)$.

\begin{theoremstar}
The number $\# X_\lambda(\mathbb F_q)$ of rational points is independent of
$\lambda$ (for $\mathbb F_q$ containing fields of definition of $\mathcal X$, $\pi$, $\lambda$
and a finite number of auxilliary varieities).
\end{theoremstar}

\begin{proof}
First note that $X_\lambda$ is isomorphic to $X_{t\lambda}$ for $t\in k^*$. 
Therefore, it is enough to show that $\# X_0(\mathbb F_q)$ is
equal to $\# X_1(\mathbb F_q)$.

Let $\bigsqcup \mathcal F_\alpha$ be the decomposition of the fixed
point set $\mathcal X^{\mathbb G_m}$ into connected components. Each
$\mathcal F_\alpha$ is a nonsingular projective variety. Moreover, $\mathcal F_\alpha$
is contained in $X_0$. (We have used the assumption (1).)

We consider the Bialynicki-Birula decomposition of $\mathcal X$ with 
respect to the $\mathbb G_m$-action:
\begin{equation*}
   \mathcal X = \bigsqcup_\alpha \mathcal X_\alpha,
\end{equation*}
where $\mathcal X_\alpha = \{ x\in \mathcal X \mid \lim_{t\to 0}
t\cdot x\in \mathcal F_\alpha\}$. By the assumption (2), the right
hand side coincides with the whole space $\mathcal X$.
It is known that the natural projection $\mathcal X_\alpha\to \mathcal
F_\alpha$ is an affine fibration whose fiber is isomoprhic to the
direct sum of positive weight space in the tangent space at $\mathcal
F_\alpha$. Therefore, we have
\begin{equation*}
   \# \mathcal X(\mathbb F_q)
   = \sum_\alpha \#\mathcal X_\alpha(\mathbb F_q)
   = \sum_\alpha \#\mathcal F_\alpha(\mathbb F_q) q^{n_\alpha},
\end{equation*}
where $n_\alpha$ is the dimension of the fiber.

We also consider the Bialynicki-Birula decomposition of $X_0$:
\begin{equation*}
   X_0 = \bigsqcup_\alpha (X_0)_\alpha.
\end{equation*}
Then $(X_0)_\alpha$ is also an affine fibration over the {\it same\/}
base $\mathcal F_\alpha$. The tangent space of $\mathcal X$ (at a point 
in $\mathcal F_\alpha$) decompose into the sum of the tangent space of 
$X_0$ (fiber direction) and $\mathbb A$ (base direction). Therefore,
the dimension of the fiber is equal to $n_\alpha - 1$. Thus
\begin{equation*}
   \# X_0(\mathbb F_q)
   = \sum_\alpha \# (X_0)_\alpha(\mathbb F_q)
   = \sum_\alpha \#\mathcal F_\alpha(\mathbb F_q) q^{n_\alpha-1}
   = \frac 1q \# \mathcal X(\mathbb F_q).
\end{equation*}
On the other hand,
\begin{equation*}
   \# \mathcal X(\mathbb F_q) = 
   \sum_{\lambda\in\mathbb F_q} \# X_\lambda(\mathbb F_q)
   = (q-1) \# X_1(\mathbb F_q) + \# X_0(\mathbb F_q).
\end{equation*}
Therefore the conclusion follows.
\end{proof}

\ifx\undefined\bysame
\newcommand{\bysame}{\leavevmode\hbox to3em{\hrulefill}\,}
\fi

\end{document}